\documentclass{asl} 

\usepackage{mathrsfs}
\usepackage[all]{xy}

\title{``Iff'' is not expressible in independence-friendly logic}

\author{Allen L. Mann}
\revauthor{Mann, Allen L.}

\address{Department of Mathematics \\
Colgate University \\
13 Oak Drive \\
Hamilton, NY 13346, USA}

\email{allen.l.mann@gmail.edu}

\genaddr{Correspondence address} 
{Allen L. Mann\\ 
407 E. Spruce St.\\ 
Caldwell, ID 83605, USA}

\newtheorem{thm}{Theorem}[section] 
\newtheorem{prop}[thm]{Proposition} 
\newtheorem{lem}[thm]{Lemma}

\theoremstyle{definition} 
\newtheorem{defn}[thm]{Definition}

\input{mann.mac}

\begin{document}

\begin{abstract}
Ordinary first-order logic has the property that two formulas $\phi$ and $\psi$ have the same meaning in a structure if and only if the formula $\phi \iff \psi$ is true in the structure. We prove that independence-friendly logic does not have this property.
\end{abstract}

\maketitle

\section{Introduction}

The \emph{meaning} of a first-order formula $\phi$ in a structure $\A$ is just the set of valuations that make the formula true in $\A$. That is,
\[
\phi^\A = \setof{\vec a \in \^NA}{\A \models \phi[\vec a]},
\]
where $A$ is the universe of $\A$, and $N$ is the number of variables in $\phi$. Given a structure $\A$ and any two first-order formulas $\phi$ and $\psi$, 
\[
\phi^\A = \psi^\A 
\]
if and only if \(\A \models \phi \iff \psi\). Thus first-order logic is able to express the concept of ``if and only if.''

Independence-friendly logic (IF logic) is a conservative extension of first-order logic that has the same expressive power as existential second-order logic \cite{Hintikka:1989, Hintikka:1996}. In IF logic the truth of a sentence is defined via a game between two players, \abelard\ ($\forall$) and \eloise\ ($\exists$). The additional expressivity is obtained by modifying the quantifiers and connectives of an ordinary first-order sentence in order to restrict the information available to the existential player, \eloise, in the associated semantic game.

In IF logic, only the information available to \eloise\ is restricted, which means existential quantifiers are not dual to universal quantifiers. To compensate, negation symbols are only allowed before atomic formulas. Generalized independence-friendly logic (IFG logic) is a variant of independence-friendly logic in which the information available to both players can be restricted, making existential quantifiers dual to universal quantifiers and allowing any formula to be negated \cite{Dechesne:2005}.

Since there are IFG-sentences that are neither true nor false, it is unclear whether IFG logic can express the concept of ``if and only if.'' For instance, one can define \(\phi \hiff{J} \psi\) as an abbreviation for the formula
\[
({\hneg\phi} \hor{J} \psi) \hand{J} (\phi \hor{J} {\hneg\psi})
\]
(the subscripts indicate what information is unavailable to the players at each move), but does this formula assert that $\phi$ and $\psi$ are logically equivalent? Is \(\phi \hiff{J} \psi\) true in a structure exactly when $\phi$ and $\psi$ have the same meaning in that structure? If not, is there some other syntactical combination of $\phi$ and $\psi$ that is true exactly when $\phi$ and $\psi$ have the same meaning? The answer to all of these questions is no.


\section{IFG logic}

\begin{defn}
Given a first-order signature $\sigma$, an \emph{atomic IFG-formula}\index{IFG-formula!atomic|mainidx} is a pair $\tuple{\phi, X}$ where $\phi$ is an atomic first-order formula and $X$ is a finite set of variables that includes every variable that appears in $\phi$ (and possibly more).
\end{defn}

\begin{defn}
Given a first-order signature $\sigma$, the language $\mathscr L_\mathrm{IFG}^\sigma$\index{$\mathscr L_\mathrm{IFG}^\sigma$|mainidx}\index{IFG-formula|(} is the smallest set of formulas such that:
\begin{enumerate}
	\item Every atomic IFG-formula is in $\mathscr L_\mathrm{IFG}^\sigma$.
	\item If $\tuple{\phi, Y}$ is in $\mathscr L_\mathrm{IFG}^\sigma$ and \(Y \subseteq X\), then $\tuple{\phi, X}$ is in $\mathscr L_\mathrm{IFG}^\sigma$.
	\item If $\tuple{\phi, X}$ is in $\mathscr L_\mathrm{IFG}^\sigma$, then $\tuple{\hneg\phi, X}$ is in $\mathscr L_\mathrm{IFG}^\sigma$.
	\item If $\tuple{\phi, X}$ and $\tuple{\psi, X}$ are in $\mathscr L_\mathrm{IFG}^\sigma$, and \(Y \subseteq X\), then $\tuple{\phi \hor{Y} \psi, X}$ is in $\mathscr L_\mathrm{IFG}^\sigma$.
	\item If $\tuple{\phi, X}$ is in $\mathscr L_\mathrm{IFG}^\sigma$, \(x \in X\), and \(Y \subseteq X\), then $\tuple{\hexists{x}{Y}\phi, X}$ is in $\mathscr L_\mathrm{IFG}^\sigma$.
\end{enumerate}
Above $X$ and $Y$ are finite sets of variables.
\end{defn}

From now on we will make certain assumptions about IFG-formulas that will allow us to simplify our notation. First, we will assume that the set of variables of $\mathscr L_\mathrm{IFG}^\sigma$ is $\setof{v_n}{n \in \omega}$. Second, since it does not matter much which particular variables appear in a formula, we will assume that variables with smaller indices are used before variables with larger indices. More precisely, if $\tuple{\phi, X}$ is a formula, \(v_j \in X\), and \(i \leq j\), then \(v_i \in X\). By abuse of notation, if $\tuple{\phi, X}$ is a formula and \(\abs X = N\), then we will say that $\phi$ has $N$ variables and write $\phi$ for $\tuple{\phi, X}$. As a shorthand, we will call $\phi$ an IFG$_N$-formula\index{IFG$_N$-formula}. Let 
\[
\mathscr L^\sigma_{\mathrm{IFG}_N} = \setof{\phi \in \mathscr L^\sigma_{\mathrm{IFG}}}{\phi \text{ has $N$ variables}}\index{$\mathscr L^\sigma_{\mathrm{IFG}_N}$}.
\] 
Third, sometimes we will write $\phi \hor{J} \psi$ instead of $\phi \hor{Y} \psi$ and $\hexists{v_n}{J}\phi$ instead of $\hexists{v_n}{Y}\phi$, where \(J = \setof{j}{v_j \in Y}\). Finally, we will use \(\phi \hand{J} \psi\) to abbreviate \(\hneg(\hneg\phi \hor{J} \hneg\psi)\) and \(\hforall{v_n}{J} \phi\) to abbreviate \(\hneg\hexists{v_n}{J}\hneg\phi\).



Truth and falsity for IFG-sentences are defined in terms of a two-player, win-loss game of imperfect information. Given an IFG-sentence, \eloise's goal is to verify the sentence, and \abelard's goal is to falsify it. The sentence is true if \eloise\ has a winning strategy, and it is false if \abelard\ has a winning strategy. For example, consider a structure $\A$ with universe $A$ and the ordinary first-order sentence
\[
\forall v_0 \exists v_1[v_0 \not= v_1].
\]
First \abelard\ chooses an element $a$ to be the value of the variable $v_0$, then \eloise\ chooses an element $b$ to be the value of the variable $v_1$. If \(a \not= b\), \eloise\ wins; otherwise, \abelard\ wins. If $\A$ has more than one element \eloise\ can win every play of the game; hence the sentence is true. If $\A$ has only one element, then \abelard\ will win every play; hence the sentence is false. Now consider the IFG$_2$-sentence 
\[
\hforall{v_0}{\emptyset}\hexists{v_1}{\set{v_0}}[v_0 \not= v_1].
\]
The subscripts indicate what information is unavailable to the players at each move. The game begins as before with \abelard\ choosing an \(a \in A\) to be the value of $v_0$. Next \eloise\ chooses an element \(b \in A\) to be the value of $v_1$, but this time she must make her choice in ignorance of the value of $v_0$. Let us assume $A$ has more than one element. On the one hand, \eloise\ does not have a winning strategy because she might blindly choose the same element as \abelard. Therefore the sentence is not true. On the other hand, \abelard\ does not have a winning strategy either, because \eloise\ might get lucky and choose a different element than the one he chose. Therefore the sentence is not false. 

Now we define the game semantics for formulas with free variables. Consider the IFG$_2$-formula 
\[
\hexists{v_1}{\set{v_0}}[v_0 \not= v_1].
\]
In order for us to decide who wins a given play of the semantic game, at the end of the game every variable must have a value. Since $v_0$ is free, neither player has the opportunity to choose its value. To get around this difficulty, before the game begins we will assign random values to all the free variables. In fact, instead of assigning values only to the free variables, we will assign values to all of the variables. Thus the first move of the game is to choose a valuation \(\vec a \in \^NA\). Play proceeds with the players modifying the initial valuation until an atomic formula is reached, at which point the the game ends and the final valuation is used to determine the winner. In the above example, play begins with values for $v_0$ and $v_1$ being chosen at random. Then \eloise\ attempts to modify the value of $v_1$ so that it is different from the value of $v_0$. Unfortunately for her, she is not allowed to see the value of $v_0$, so her task is no easier than before. However, suppose an oracle revealed to \eloise\ that the initial valuation belonged to a subset $V$ of the space of all valuations $\^2A$. \eloise\ might be able to use this information to devise a winning strategy. For example, suppose \(A = \set{0,1}\). Then \(\^2A = \set{00, 01, 10, 11}\). If the oracle tells \eloise\ that the initial valuation belongs to the set \(V = \set{00, 01}\), then \eloise\ will know to choose 1 for the value of $v_1$. Thus \eloise\ has a winning strategy for the game that begins by choosing the initial valuation from $V$ instead of from $\^2A$. A set of valuations, such as $V$, is called a \emph{team}. We say that the formula \(\hexists{v_1}{\set{v_0}}[v_0 \not= v_1]\) is true in $\A$ relative to $V$, and that $V$ is a winning team for $\phi$ in $\A$.

Disjunctions and conjunctions are moves for the players, as well. In the game corresponding to the formula \(\psi_1 \hor{Y} \psi_2\), \eloise\ must choose which disjunct she wishes to verify without knowing the values of the variables in $Y$. Dually, in the game corresponding to \(\psi_1 \hand{Y} \psi_2\), \abelard\ chooses which conjunct \eloise\ must verify, but his choice is not allowed to depend on the variables in $Y$.

Negation is handled by having the players switch roles. \eloise\ attempts to verify ${\hneg\psi}$ by falsifying $\psi$, and \abelard\ attempts to falsify ${\hneg\psi}$ by verifying $\psi$.

In general, if $\phi$ is an IFG$_N$-formula and \(V,W \subseteq \^NA\) are teams, then $\phi$ is true in $\A$ relative to $V$, denoted \(\A \modelt \phi[V]\), if and only if \eloise\ has a winning strategy for the semantic game, given she knows the initial valuation belongs to $V$. Dually, $\phi$ is false in $\A$ relative to $W$, denoted \(\A \modelf \phi[W]\), if and only if \abelard\ has a winning strategy, given he knows the initial valuation belongs to $W$. In the first case, we say that $V$ is a \emph{winning team} (or \emph{trump}) for $\phi$ in $\A$, and in the second case, that $W$ is a \emph{losing team} (or \emph{cotrump}) for $\phi$ in $\A$.

Finally, we need to connect the game semantics for IFG-sentences with the game semantics for IFG-formulas. If $\phi$ is an IFG-sentence, then the initial valuation is irrelevant because the value of every variable is modified during the course of the game. If $\phi$ has free variables, then in order to have a winning strategy, a player must be able to win no matter what the initial valuation is. Therefore we define \(\A \modelt \phi\) if and only if \(\A \modelt \phi[\^NA]\) and \(\A \modelf \phi\) if and only if \(\A \modelf \phi[\^NA]\). In the future, we will abbreviate similar statements by writing \(\A \modeltf \phi\) if and only if \(\A \modeltf \phi[\^NA]\).

It is worth noting that restricting the information available to the players does not affect their moves, only their strategies. Therefore, restricting the information available to one player does not help his or her opponent. \eloise\ has a winning strategy if and only if she wins regardless of how \abelard\ plays. Withholding information from her makes it harder for her to have a winning strategy, but withholding information from \abelard\ does not make it easier.

We hope this summary of the game semantics for IFG logic is sufficient. A more rigorous treatment can be found in \cite[Section 1.2]{Mann:2007a} or \cite[Section 1.3]{Mann:2007}.

\section{Trump semantics}

Wilfrid Hodges made an important breakthrough when he found a way to define a Tarski-style semantics for independence-friendly logic \cite{Hodges:1997a, Hodges:1997b}. We now recall the necessary details.

\begin{defn} 
Two valuations \(\vec a, \vec b \in\, \^NA\) \emph{agree outside of \(J \subseteq N\)}, denoted \(\vec a \approx_J \vec b\)\index{\(\vec a \approx_J \vec b\)\quad $\vec a$ and $\vec b$ agree outside of $J$|mainidx}, if 
\[
\vec a \restrict (N\setminus J) = \vec b \restrict(N\setminus J).
\]
\end{defn}

\begin{defn} 
Given any set $V$, a \emph{cover} of $V$ is a collection of sets $\mathscr U$ such that \(V = \bigcup\mathscr U\). A \emph{disjoint cover} is a cover whose members are pairwise disjoint.
\end{defn}

\begin{defn} 
Let \(V \subseteq \^NA\), and let $\mathscr U$ be a cover of $V$. Then $\mathscr U$ is called \emph{$J$-saturated}\index{cover!$J$-saturated|mainidx} if every \(U \in \mathscr U\) is closed under $\approx_J$. That is, for every \(\vec a, \vec b \in V\), if \(\vec a \approx_J \vec b\) and \(\vec a \in U \in \mathscr U\), then \(\vec b \in U\).
\end{defn}

\begin{defn} 
Define a partial operation $\bigcup_J$\index{$\bigcup_J \mathscr U$|mainidx} on sets of teams by declaring \(\bigcup_J \mathscr U = \bigcup \mathscr U\) whenever $\mathscr U$ is a $J$-saturated disjoint cover of $\bigcup \mathscr U$ and letting $\bigcup_J \mathscr U$ be undefined otherwise. Thus the formula \(V = \bigcup_J \mathscr U\) asserts that $\mathscr U$ is a $J$-saturated disjoint cover of $V$. We will use the notation \(V_1 \cup_J V_2\)\index{$V_1 \cup_J V_2$|mainidx} to abbreviate \(\bigcup_J \set{V_1, V_2}\).
\end{defn}

\begin{defn} 
A function \(f\colon V \to A\) is \emph{independent of $J$}, denoted \(f\colon V \toind{J} A\)\index{$f\colon V \toind{J} A$ \quad $f$ is independent of $J$|mainidx},  if \(f(\vec{a}) = f(\vec{b})\) whenever \(\vec{a} \approx_J \vec{b}\). 
\end{defn}

\begin{defn} 
Let \(\vec a \in \^NA\). For every \(n < N\) and \(b \in A\), define $\vec a(n:b)$ to be the valuation that is like $\vec a$ except the $n$th value has been changed to $b$, i.e.,
\[
\vec a(n:b) = \vec a \restrict(N\setminus\set{n}) \cup \set{\tuple{n,b}}.
\]
Let \(V,W \subseteq \^NA\) and \(f\colon V \to A\). Define
\begin{align*}
V(n:f) &= \setof{\vec a(n:f(a))}{\vec a \in V}, \\
W(n:A) &= \setof{\vec a(n:b)}{\vec a \in W,\, b \in A}.
\end{align*}

\end{defn}

The next theorem has appeared in many different forms in the literature.
Hodges' original formulation (for IF logic) appears in \cite[Theorem 7.5]{Hodges:1997a}. Dechesne's version for IFG logic appears in \cite[Theorem 5.3.5]{Dechesne:2005}.

\begin{thm}[Hodges] 
\label{trump semantics}
\index{$\A \modeltf \phi[V]$|mainidx}
Let $\phi$ be an IFG$_N$-formula, let $\A$ be a suitable structure, and let \(V,W \subseteq \^NA\).
\begin{itemize}
	\item If $\phi$ is atomic, then  
	\begin{itemize}
		\item[(+)] \(\A \modelt \phi[V]\) if and only if for every \(\vec a \in V\), \(\A \models \phi[\vec a]\),
		\item[($-$)] \(\A \modelf \phi[W]\) if and only if for every \(\vec b \in W\), \(\A \not\models \phi[\vec b]\).	
	\end{itemize}
	\item If $\phi$ is ${\hneg\psi}$, then 
	\begin{itemize}
		\item[(+)] \(\A \modelt {\hneg\psi[V]}\) if and only if \(\A \modelf \psi[V]\),
		\item[($-$)] \(\A \modelf {\hneg\psi[W]}\) if and only if \(\A \modelt \psi[W]\).
	\end{itemize}
	\item If $\phi$ is $\psi_1 \hor{J} \psi_2$, then 
	\begin{itemize}
		\item[(+)] \(\A \modelt \psi_1 \hor{J} \psi_2[V]\) if and only if \(\A \modelt \psi_1[V_1]\) and \( \A \modelt \psi_2[V_2]\) for some \(V = V_1 \cup_J V_2\),
		\item[($-$)] \(\A \modelf \psi_1 \hor{J} \psi_2[W]\) if and only if \(\A \modelf \psi_1[W]\) and \( \A \modelf \psi_2[W]\).
	\end{itemize}
	\item If $\phi$ is $\hexists{v_n}{J}\psi$, then 
	\begin{itemize}
		\item[(+)] \(\A \modelt \hexists{v_n}{J}\psi[V]\) if and only if \(\A \modelt \psi[V(n:f)]\) for some \(f\colon V \toind{J} A\),
		\item[($-$)] \(\A \modelf \hexists{v_n}{J}\psi[W]\) if and only if \(\A \modelf \psi[W(n:A)]\).
	\end{itemize}
\end{itemize}
\end{thm}

\begin{proof}
By two simultaneous inductions on the subformulas of $\phi$. A full proof using the present notation can be found in \cite[Theorem 1.32]{Mann:2007}.
\end{proof}





\section{IFG-cylindric set algebras}

We introduced IFG-cylindric set algebras in \cite{Mann:2007a, Mann:2007} as a way to study the algebra of IFG logic.

Recall from \cite[p.~2]{Henkin:1981} or from \cite[Definition 4.3.4 on p.~154]{Henkin:1985} that the universe of $\Cyls_{N}(\A)$, the $N$-dimensional cylindric set algebra over $\A$, consists of the meanings of all the $N$-variable, first-order formulas expressible in the language of $\A$, where the meaning of a formula is defined by 
\[
\phi^\A = \setof{\vec a \in \^NA}{\A \models \phi[\vec a]}.
\]
Similarly, the universe of \(\Cyls_{\mathrm{IFG}_{N}}(\A)\), the $N$-dimensional IFG-cylindric set algebra over $\A$, consists of the meanings of all the IFG$_N$-formulas expressible in the language of $\A$, where the meaning of an IFG$_N$-formula is given by
\[
\norm{\phi}_\A = \tuple{\trump{\phi}_\A, \cotrump{\phi}_\A},
\]
\[
\trump{\phi}_\A = \setof{V \subseteq \^NA}{\A \modelt \phi[V]}, 
\qquad
\cotrump{\phi}_\A = \setof{W \subseteq \^NA}{\A \modelf \phi[W]}.
\]

More generally, we can define IFG-cylindric set algebras  without reference to a base structure $\A$.

\begin{defn} 
An \emph{IFG-cylindric power set algebra}\index{independence-friendly cylindric power set algebra|mainidx} is an algebra whose universe is \(\powerset(\powerset(\^NA)) \times \powerset(\powerset(\^NA))\), where $A$ is a set and $N$ is a natural number. The set $A$ is called the \emph{base set}\index{base set|mainidx}, and the number $N$ is called the \emph{dimension}\index{dimension|mainidx} of the algebra. Every element $X$ of an IFG-cylindric power set algebra is an ordered pair of sets of teams. We will use the notation $X^+$\index{$X^+$ \quad truth coordinate of $X$|mainidx} to refer to the first coordinate of the pair, and $X^-$\index{$X^-$ \quad falsity coordinate of $X$|mainidx} to refer to the second coordinate. There are a finite number of operations: 

\begin{itemize}
	\item the constant \(0 = \tuple{\set{\emptyset}, \powerset(\^NA)}\);
	\item the constant \(1 = \tuple{\powerset(\^NA), \set{\emptyset}}\);
	\item for all \(i,j < N\), the constant $D_{ij}$\index{$D_{ij}$ \quad diagonal element|mainidx} is defined by 
	\begin{itemize}
		\item[(+)] \(D_{ij}^+ = \powerset(\setof{\vec a \in\, \^NA}{a_i = a_j})\),
		\item[($-$)] \(D_{ij}^- = \powerset(\setof{\vec a \in\, \^NA}{a_i \not= a_j})\);
	\end{itemize}
	\item if \(X = \tuple{X^+,\, X^-}\), then \(\n{X}\index{$\n{X}$ \quad negation of $X$|mainidx} = \tuple{X^-, X^+}\);
	\item for every \(J \subseteq N\), the binary operation $+_J$\index{$X +_J Y$|mainidx} is defined by
	\begin{itemize}
		\item[(+)] \(V \in (X +_J Y)^+\) if and only if \(V = V_1 \cup_J V_2\) for some
				 \(V_1 \in X^+\) and \(V_2 \in Y^+\),
		\item[($-$)] \((X +_J Y)^- = X^- \cap Y^-\);
	\end{itemize}
	\item for every \(J \subseteq N\), the binary operation $\cdot_J$\index{$X \cdot_J Y$|mainidx} is defined by
	\begin{itemize}
		\item[(+)] \((X \cdot_J Y)^+ = X^+ \cap Y^+\),
		\item[($-$)] \(W \in (X \cdot_J Y)^-\) if and only if \(W = W_1 \cup_J W_2\) for some
				 \(W_1 \in X^-\) and \(W_2 \in Y^-\);
	\end{itemize}
	\item for every \(n < N\) and \(J \subseteq N\), the unary operation $C_{n,J}$\index{$C_{n,J}(X)$ \quad cylindrification of $X$|mainidx} is defined by
		\begin{itemize}
			\item[(+)] \(V \in C_{n,J}(X)^+\) if and only if \(V(n:f) \in X^+\) for some 
					\(f\colon V \toind{J} A\),
			\item[($-$)] \(W \in C_{n,J}(X)^- \) if and only if \(W(n:A) \in X^-\).
		\end{itemize}
\end{itemize}
\end{defn}

\begin{defn} 
An \emph{IFG-cylindric set algebra}\index{independence-friendly cylindric set algebra|mainidx} (or \emph{IFG-algebra}, for short) is any subalgebra of an IFG-cylindric power set algebra. An \emph{IFG$_N$-cylindric set algebra}\index{IFG$_N$-cylindric set algebra|mainidx} (or \emph{IFG$_N$-algebra}) is an IFG-cylindric set algebra of dimension $N$.
\end{defn}

The operations $+_\emptyset$ and $+_N$ are of particular interest. Since every disjoint cover of $V$ is $\emptyset$-saturated, \(V \in (X +_\emptyset Y)^+\) if and only if there is a disjoint cover \(V = V_1 \cup V_2\) such that \(V_1 \in X^+\) and \(V_2 \in Y^+\). At the other extreme, \(V = V_1 \cup_N V_2\) if and only if \(V_1 = V\) and \(V_2 = \emptyset\) or vice versa.

Also, the element \(\Omega = \tuple{\set{\emptyset}, \set{\emptyset}}\) is present in most, but not all, IFG-algebras.

\section{Suits and double suits}

Meanings of IFG-formulas have the property that \(\trump{\phi} \cap\ \cotrump{\phi} = \set{\emptyset}\), and \(V' \subseteq V \in \norm{\phi}^\pm\) implies \(V' \in \norm{\phi}^\pm\). These facts inspire the following definitions.

\begin{defn} 
A nonempty set $X^* \subseteq\powerset(\^NA)$ is called a \emph{suit}\index{suit|mainidx} if \(V' \subseteq V \in X^*\) implies \(V' \in X^*\). A \emph{double suit}\index{double suit|mainidx} is a pair $\tuple{X^+, X^-}$ of suits such that \(X^+ \cap X^- = \set{\emptyset}\).
\end{defn}

\begin{defn} 
An IFG-algebra is \emph{suited}\index{IFG-cylindric set algebra!suited|mainidx} if all of its elements are pairs of suits. It is \emph{double-suited}\index{IFG-cylindric set algebra!double-suited|mainidx} if all of its elements are double suits.
\end{defn}

\begin{prop}[Proposition 2.10 in \cite{Mann:2007a}] 
\label{subalgebra generated suited}
The IFG$_N$-algebra generated by a set of pairs of suits is a suited IFG$_N$-algebra.
\end{prop}

\begin{prop}[Proposition 2.11 in \cite{Mann:2007a}] 
\label{subalgebra generated double-suited}
The IFG$_N$-algebra generated by a set of double suits is a double-suited IFG$_N$-algebra. In particular, \(\Cyls_{\mathrm{IFG}_{N}}(\A)\) is a double-suited IFG$_N$-algebra.
\end{prop}

For the rest of the paper we will only be concerned with double suits. The next proposition is a summary of results from \cite[Section 2.5]{Mann:2007a}.

\begin{prop}\label{section 2.5}
If $X$ and $Y$ are double suits,
\begin{enumerate}
	\item \(X +_J 0 = X = X \cdot_J 1\),
	\item \(X +_J 1 = 1\) and \(X \cdot_J 0 = 0\),
	\item \(X +_N Y = \tuple{X^+ \cup Y^+,\, X^- \cap Y^-}\) and \(X \cdot_N Y = \tuple{X^+ \cap Y^+,\, X^- \cup Y^-}\).
\end{enumerate}
\end{prop}

Part (c) implies that $+_N$ and $\cdot_N$ are lattice operations. Thus we can define a partial order on any double-suited IFG$_N$-algebra by declaring \(X \leq Y\) if and only if \(X +_N Y = Y\). It follows that \(X \leq Y\) if and only if \(X^+ \subseteq Y^+\) and \(Y^- \subseteq X^-\).

\begin{prop}\label{X < Omega,Y}
If $X$ and $Y$ are double suits, \(X \leq \Omega\), and \(X \leq Y\), then \(X +_J Y = Y\).
\end{prop}

\begin{proof}
If \(X \leq \Omega\), then \((X +_J Y)^+ = Y^+\). To see why, suppose \(X \leq \Omega\) and \(V \in (X +_J Y)^+\). Then \(V = V_1 \cup_J V_2\) for some \(V_1 \in X^+\) and \(V_2 \in Y^+\), but since \(X^+ = \set{\emptyset}\) we must have \(V_1 = \emptyset\) and \(V_2 = V\). Hence \(V \in Y^+\). Conversely, suppose \(V \in Y^+\). Then \(V = \emptyset \cup_J V\), where \(\emptyset \in X^+\) and \(V \in Y^+\), so \(V \in (X +_J Y)^+\).

If \(X \leq Y\), then \(Y^- \subseteq X^-\), so \((X +_J Y)^- = X^- \cap Y^- = Y^-\).
\end{proof}

Whereas an ordinary cylindric algebra is an expansion of a Boolean algebra, we should not expect the same to be true for IFG-algebras because of the failure of the law of excluded middle in IFG logic. Somewhat miraculously, double-suited IFG-algebras do have an underlying structure that is as close to being a Boolean algebra as possible without satisfying the complementation axioms.

\begin{defn} 
A \emph{De Morgan algebra}\index{De Morgan algebra|mainidx} \(\A = \seq{A; 0, 1, {\hneg}\,, \join, \meet}\) is a bounded distributive lattice with an additional unary operation $\hneg\ $ that satisfies \(\hneg{\hneg x} = x\) and \(\hneg(x \join y) = {\hneg x} \meet {\hneg y}\).
\end{defn}

\begin{defn} 
A \emph{Kleene algebra}\index{Kleene algebra|mainidx} is a De Morgan algebra that satisfies the additional axiom \(x \meet {\hneg x} \leq y \join {\hneg y}\).
\end{defn}

\begin{thm}[Theorem 2.31 in \cite{Mann:2007a}] 
\label{Kleene reduct}
The reduct of a double-suited IFG$_N$-algebra to the signature \(\tuple{0, 1, \n{ }, +_N, \cdot_N}\) is a Kleene algebra.
\end{thm} 

Given a set $A$, let $\Suit(\^NA)$\index{$\Suit(\^NA)$ \quad set of suits in \(\powerset(\powerset(\^NA))\)} denote the set of all suits in \(\powerset(\powerset(\^NA))\), and let $\DSuit(\^NA)$\index{$\DSuit(\^NA)$ \quad set of double-suits in \(\powerset(\powerset(\^NA)) \times \powerset(\powerset(\^NA))\)} denote the set of all double suits in \(\powerset(\powerset(\^NA)) \times \powerset(\powerset(\^NA))\). Since the meaning of every IFG-formula is a double suit, the universe of \(\Cyls_{\mathrm{IFG}_{N}}(\A)\) is contained in \(\DSuit(\^NA)\). Therefore \(\abs{\DSuit(\^NA)}\) gives an upper bound for the size of $\Cyls_{\mathrm{IFG}_{N}}(\A)$. Cameron\index{Cameron, P.} and Hodges\index{Hodges, W.} \cite{Cameron:2001} count suits and double suits in the case when \(N = 1\). The results of their calculations are shown in Table \ref{counting suits and double suits}, where \(m = \abs{A}\), \(f(m) = \abs{\Suit(A)}\), and \(g(m) = \abs{\DSuit(A)}\). They remark that ``one can think of the ratio of $g(m)$ to $2^m$ as measuring the expressive strength of [IFG logic] compared with ordinary first-order logic---always bearing in mind that [IFG logic] may have a rather unorthodox notion of what is worth expressing'' \cite[p.~679]{Cameron:2001}. Cameron\index{Cameron, P.} and Hodges\index{Hodges, W.} also prove that given any finite set $A$, there is a structure $\A$ such that the universe of $\Cyls_{\mathrm{IFG}_{N}}(\A)$ is exactly $\DSuit(\^NA)$ \cite[Corollary 3.4]{Cameron:2001}.

\begin{table} 
\[
\left.
\begin{array}{c|l|l|l}
m & 2^m		& f(m)												& g(m) \\\hline 
0 & 1		& 1													& 1 \\
1 &	2		& 2													& 3 \\
2 & 4		& 5													& 11 \\
3 & 8		& 19												& 55 \\
4 & 16		& 167												& 489 \\
5 & 32		& 7,\!580											& 17,\!279 \\
6 & 64		& 7,\!828,\!353										& 15,\!758,\!603 \\
7 & 128		& 2,\!414,\!682,\!040,\!997						& 4,\!829,\!474,\!397,\!415 \\
8 & 256		& 56,\!130,\!437,\!228,\!687,\!557,\!907,\!787	& 112,\!260,\!874,\!496,\!010,\!913,\!723,\!317
\end{array}
\right.
\]
  \caption{Counting suits and double suits}
  \label{counting suits and double suits}
\end{table}

\begin{prop} 
\label{all double suits}
Let $\A$ be a finite structure with at least two elements, and in which every element is named by a constant symbol. Then the universe of \(\Cyls_{\mathrm{IFG}_{N}}(\A)\) is exactly \(\DSuit(\^NA)\).\index{$\DSuit(\^NA)$ \quad double-suited IFG$_N$-cylindric set algebra over $A$}
\end{prop}

\begin{proof}
For every \(\vec a \in \^NA\), let $\phi_{\vec a}$ be the formula \(v_0 = a_0 \hand{\emptyset} \cdots \hand{\emptyset} v_{N-1} = a_{N-1}\). For every \(V \subseteq \^NA\), let $\phi_V$ be the formula \(\bigvee_{\!/\emptyset} \setof{\phi_{\vec a}}{\vec a \in V}\). 
Then
\begin{align*}
\norm{\phi_{\vec a}} &= \tuple{\powerset(\set{\vec a}), \powerset(\^NA \setminus \set{\vec a})}, \\
\norm{\phi_V} &= \tuple{\powerset(V), \powerset(\^NA \setminus V)}.
\end{align*}

Let $X$ be a double suit, and let \(X^+ = \powerset(V_0) \cup \cdots \cup \powerset(V_{k-1})\). Let $\phi$ be the formula 
\(
\phi_{V_0} \hor{N} \cdots \hor{N} \phi_{V_{k-1}}.
\)
Then \(\trump{\phi} = X^+\). Similarly, there is a formula $\psi$ such that \(\cotrump{\psi} = X^-\). Let $c$ be a constant symbol naming one of the elements of $\A$, let $\chi$ be the formula \(v_0 = c \hor{N} v_0 \not= c\), and let \(V = V_0 \cup \cdots \cup V_{k-1}\). Then \(\norm{\chi} = \Omega\), and 
\begin{align*}
X' &= \norm{\phi \hor{N} \chi} = \tuple{X^+, \set{\emptyset}}, \\
X'' &= \norm{\psi \hand{N} \chi} = \tuple{\set{\emptyset}, X^-}, \\
Y &= \norm{\phi_V} = \tuple{\powerset(V), \powerset(\^NA \setminus V)}.
\end{align*}
It suffices to show that
\begin{align*}
\norm{(\phi \hor{N} \chi) \hand{N} ((\psi \hand{N} \chi) \hor{N} \phi_V)} 
&= X' \cdot_N (X'' +_N Y) \\
&= \tuple{X^+ \cap \powerset(V),\ X^- \cap \powerset(\^NA \setminus V)} \\
&= \tuple{X^+\!,\, X^-} \\
&= X.
\end{align*}
Therefore \(X \in \Cyls_{\mathrm{IFG}_{N}}(\A)\).
\end{proof}

At this point, it is natural to ask which double suits can be the meanings of ordinary first-order formulas (that is, IFG-formulas whose independence sets are all empty). Ordinary first-order formulas have the property that \(\A \modelt \phi[V]\) if and only if \(\A \modelt \phi[\set{\vec a}]\) for all \(\vec a \in V\). Hence if \(\A \modelt \phi[V]\) and \(\A \modelt \phi[V']\), then \(\A\modelt \phi[V \cup V']\). It follows that the set of winning teams for an ordinary first-order formula $\phi$ is simply the power set of the set of valuations that satisfy $\phi$. That is,
\[
\trump{\phi}_\A = \powerset(\phi^\A).
\]
Ordinary first-order formulas also have the property that for every \(\vec a \in \^NA\) either \(\A \modelt \phi[\set{\vec a}]\) or \(\A \modelf \phi[\set{\vec a}]\). These facts inspire the following definitions.

\begin{defn} 
A double suit $X$ is \emph{flat}\index{suit!flat|mainidx} if there is a \(V \subseteq \^NA\) such that \(X^+ = \powerset(V)\).
\end{defn}

\begin{prop}\label{flat absorption}
If $X$ and $Y$ are double suits, \(X \leq Y\), and $Y$ is flat, then \(X +_J Y = Y\).
\end{prop}

\begin{proof}
Suppose \(X \leq Y\) and \(Y^+ = \powerset(V)\). If \(V' \in (X +_J Y)^+\), then \(V' = V_1 \cup_J V_2\) for some \(V_1 \in X^+\) and \(V_2 \in Y^+\). Hence \(V_1, V_2 \subseteq V\) because \(X^+ \subseteq Y^+ = \powerset(V)\). Thus \(V' \subseteq V\), which implies \(V' \in Y^+\). Conversely, if \(V' \in Y^+\), then \(V' = \emptyset \cup_J V'\), where \(\emptyset \in X^+\) and \(V' \in Y^+\), so \(V' \in (X +_J Y)^+\). 

Also, since \(Y^- \subseteq X^-\) we have \((X +_J Y)^- = X^- \cap Y^- = Y^-\).
\end{proof}

\begin{defn} 
A double suit $X$ is \emph{perfect}\index{double suit!perfect|mainidx} if there is a \(V \subseteq \^NA\) such that 
\[
X = \tuple{\powerset(V),\, \powerset(\^NA \setminus V)}.
\]
\end{defn}

In \cite{Mann:2008aa}, we showed that an IFG-formula $\phi$ is equivalent to an ordinary first-order formula in a structure $\A$ if and only if $\norm{\phi}_\A$ is perfect. It is worth noting that $\Cyls_{\mathrm{IFG}_{N}}(\A)$ is generated by its perfect elements because it is generated by the meanings of atomic formulas, which are all perfect.

\section{$\Cyls_{\mathrm{IFG}_{1}}(\mathbf{2})$ is hereditarily simple}

Let $\mathbf{2}$ be the structure with universe $\set{0,1}$ in which both elements are named by constant symbols. Then \(\Cyls_{\mathrm{IFG}_{1}}(\mathbf{2}) = \DSuit_1(\set{0,1})\). The distributive lattice structure of $\Cyls_{\mathrm{IFG}_{1}}(\mathbf{2})$ is shown in Figure \ref{CsIFG(2)b}, where the join operation is $+_{\set{0}}$ and the meet operation is $\cdot_{\set{0}}$,
\begin{align*}
A &= \tuple{\powerset(\set{0}) \cup \powerset(\set{1}),\, \set{\emptyset}}, \\
B &= \tuple{\powerset(\set{0}),\, \set{\emptyset}}, \\
C &= \tuple{\powerset(\set{1}),\, \set{\emptyset}}, \\
\norm{v_0 = 0} &= \tuple{\powerset(\set{0}),\, \powerset(\set{1})}, \\
\norm{v_0 = 1} &= \tuple{\powerset(\set{1}),\, \powerset(\set{0})}.
\end{align*}
The goal of this section is to show that $\Cyls_{\mathrm{IFG}_{1}}(\mathbf{2})$ is hereditarily simple.

\begin{figure} 
\[
\xymatrix{
&								& 								1 \ar@{-}[d] \\
&								& 								A \ar@{-}[dl]\ar@{-}[dr] \\
&								B \ar@{-}[dl]\ar@{-}[dr]		& 									& C \ar@{-}[dl]\ar@{-}[dr] \\
\norm{v_0 = 0} \ar@{-}[dr]		&								& \Omega \ar@{-}[dl]\ar@{-}[dr]	&								& \norm{v_0 = 1} \ar@{-}[dl] \\
&								\n{C} \ar@{-}[dr]				&									& \n{B} \ar@{-}[dl] \\
&																& \n{A} \ar@{-}[d] \\
&																& 0
}
\]
  \caption{$\Cyls_{\mathrm{IFG}_{1}}(\mathbf{2})$}
  \label{CsIFG(2)b}
\end{figure}

In \cite[Proposition 2.5]{Mann:2007a}, we proved that an IFG$_N$-sentence can have one of only three possible meanings: 0, $\Omega$, and 1. Thus, if we think of $C_{0,J_0}\ldots C_{N-1,J_{N-1}}$ an a single operation that quantifies (cylindrifies) all of the variables of an IFG$_N$-formula, then the range of that operation is the IFG$_N$-algebra \(\set{0, \Omega, 1}\).

\begin{prop}[Proposition 2.51 in \cite{Mann:2007a}] 
\label{C_0,J...C_N-1,J double suit}
If $X$ is a double suit, then 
\[
C_{0,J_0}\ldots C_{N-1,J_{N-1}}(X) = 
\begin{cases}
1	    	& \text{if \(X \not\leq \Omega\)}, \\
\Omega		& \text{if \(0 < X \leq \Omega\)}, \\
0      		& \text{if \(X = 0\)}.
\end{cases}
\]
\end{prop}

\begin{lem} 
\label{congruence 0, Omega, 1}
Let $\cong$ be a congruence on a double-suited IFG$_N$-algebra. If 0, $\Omega$ (if present), or 1 are congruent to any other element, then $\cong$ is the total congruence.
\end{lem}

\begin{proof}
First we will show that if \(0 \cong 1\), then $\cong$ is the total congruence. If \(0 \cong 1\), then for every $X$ we have \(X = X +_\emptyset 0 \cong X +_\emptyset 1 = 1\). Hence $\cong$ is the total congruence. Next we will show that if \(0 \cong \Omega\) or \(1 \cong \Omega\), then $\cong$ is the total congruence. If \(0 \cong \Omega\), then \(\Omega = \n{\Omega} \cong \n{0} = 1\). Similarly, if \(1 \cong \Omega\), then \(\Omega = \n{\Omega} \cong \n{1} = 0\).

Now suppose \(0 \not= X \cong 0\). Then either 
\[
C_{0,\emptyset}\ldots C_{N-1,\emptyset}(X) = 1 
\quad \text{or} \quad 
C_{0,\emptyset}\ldots C_{N-1,\emptyset}(X) = \Omega.
\] 
In the first case, \(0 = C_{0,\emptyset}\ldots C_{N-1,\emptyset}(0) \cong C_{0,\emptyset}\ldots C_{N-1,\emptyset}(X) = 1\). In the second case, \(0 = C_{0,\emptyset}\ldots C_{N-1,\emptyset}(0) \cong C_{0,\emptyset}\ldots C_{N-1,\emptyset}(X) = \Omega\). Either way, $\cong$ is the total congruence. In addition, if \(1 \not= X \cong 1\), then \(0 \not= \n{X} \cong 0\), so $\cong$ is the total congruence.

Finally, if \(\Omega \not= X \cong \Omega\), then either \(X \not\leq \Omega\) or \(\n{X} \not\leq \Omega\). Hence either 
\[
\Omega = C_{0,\emptyset}\ldots C_{N-1,\emptyset}(\Omega) \cong C_{0,\emptyset}\ldots C_{N-1,\emptyset}(X) = 1
\]
or
\[
\Omega = C_{0,\emptyset}\ldots C_{N-1,\emptyset}(\Omega) \cong C_{0,\emptyset}\ldots C_{N-1,\emptyset}(\n{X}) = 1.
\]
Thus $\cong$ is the total congruence.
\end{proof}

\begin{lem} 
\label{congruence X < Omega < Y}
Let $\cong$ be a congruence on any double-suited IFG$_N$-algebra that includes $\Omega$. If \(X < \Omega < Y\) and \(X \cong Y\), then $\cong$ is the total congruence.
\end{lem}

\begin{proof}
If \(X < \Omega < Y\), and \(X \cong Y\), then \(\Omega = X +_N \Omega \cong Y +_N \Omega = Y\), so by the previous lemma $\cong$ is the total congruence.
\end{proof}

\begin{lem} 
\label{congruence Omega < X < Y}
Let $\cong$ be a nontrivial congruence on any double-suited IFG$_N$-algebra that includes $\Omega$. Then there exist elements $X$ and $Y$ such that \(X \cong Y\)and \(\Omega \leq X < Y\).
\end{lem}

\begin{proof}
Since $\cong$ is nontrivial, there exist distinct elements $X''$ and $Y''$ such that \(X'' \cong Y''\). Either \((X'')^+ \not= (Y'')^+\) or \((X'')^- \not= (Y'')^-\). In the first case, let \(X' = X'' +_N \Omega\) and \(Y' = Y'' +_N \Omega\); in the second case, let \(X' = \n{(X'')} +_N \Omega\) and \(Y' = \n{(Y'')} +_N \Omega\). In both cases, \((X')^+ \not= (Y')^+\), \(X' \cong Y'\), and \(\Omega \leq X',Y'\). Now let \(X = X' = X' +_N X'\) and \(Y = X' +_N Y'\). Then \(X \cong Y\) and \(\Omega \leq X < Y\).
\end{proof}

\begin{thm} 
\label{CsIFG(2) is simple}
$\Cyls_{\mathrm{IFG}_{1}}(\mathbf{2})$ is simple.
\end{thm}

\begin{proof}
By \lemref{congruence Omega < X < Y} it suffices to consider the congruences generated by pairs of elements from the interval above $\Omega$. Using the technique of perspective edges, we can see that if \(A \cong B\), then \(C \cong \Omega\), because \(A \cdot_N C = C\) and \(B \cdot_N C = \Omega\). Thus $\cong$ is the total congruence by \lemref{congruence 0, Omega, 1}. Similarly, if \(A \cong C\) then \(B \cong \Omega\). Finally, if \(B \cong C\) then \(B \cong A\) because \(B = B +_N \norm{v_0 = 0}\) and \(A = C +_N \norm{v_0 = 0}\).
\end{proof}

\begin{prop} 
\label{subalgebras of CsIFG(2)}
The proper subalgebras of $\Cyls_{\mathrm{IFG}_{1}}(\mathbf{2})$ are \(\set{0,1}\), \(\set{0, \Omega, 1}\), and those shown in Figure \ref{figure:subalgebras of CsIFG(2)}.
\begin{figure} 
\[
\xymatrix{
	1 \ar@{-}[d]		&					& 1 \ar@{-}[dl] 		&					&						& 1 \ar@{-}[dr]			&				\\
	A \ar@{-}[d] 		& 	B \ar@{-}[dr]	&						&					&						&						& C \ar@{-}[dl]	\\
	\Omega \ar@{-}[d] 	&					& \Omega \ar@{-}[dr] 	&					&						& \Omega \ar@{-}[dl]	&				\\
	\n{A} \ar@{-}[d]	&					&						& \n{B} \ar@{-}[dl]	& 	\n{C} \ar@{-}[dr] 	&						&				\\
	0					&					& 0 					&					&						& 0						&
}
\]
\[
\xymatrix{
					& 1 \ar@{-}[d]			&					&						& 1 \ar@{-}[d]			&				&						& 1 \ar@{-}[d]						&						\\
					& A \ar@{-}[dl] 		&					&						& A \ar@{-}[dr] 		&				&						& A \ar@{-}[dl] \ar@{-}[dr]			&						\\
	B \ar@{-}[dr] 	& 						&					& 						&						& C \ar@{-}[dl]	&	B \ar@{-}[dr] 		& 									& C \ar@{-}[dl]			\\
					& \Omega \ar@{-}[dr] 	&					&						& \Omega \ar@{-}[dl]	&				&						& \Omega \ar@{-}[dl] \ar@{-}[dr] 	&						\\
					&						& \n{B} \ar@{-}[dl]	& 	\n{C} \ar@{-}[dr]	&						& 				&	\n{C} \ar@{-}[dr]	&									& \n{B} \ar@{-}[dl]		\\
					& \n{A} \ar@{-}[d]		&					&						& \n{A} \ar@{-}[d]		&				&						& \n{A} \ar@{-}[d]					&						\\
					& 0						&					&						& 0						&				&						& 0									&			
}
\]
  \caption{Subalgebras of $\Cyls_{\mathrm{IFG}_{1}}(\mathbf{2})$}
  \label{figure:subalgebras of CsIFG(2)}
\end{figure}

\end{prop}

\begin{proof}
It is easy to check that the IFG$_1$-algebras \(\set{0,1}\) and \(\set{0, \Omega, 1}\) are subalgebras of \(\Cyls_{\mathrm{IFG}_{1}}(\mathbf{2})\). Consider the subalgebra \(\gen{A} = \set{0, \n{A}, \Omega, A, 1}\). Recall that since \(\Cyls_{\mathrm{IFG}_{1}}(\mathbf{2})\) is double-suited, \(X +_J 0 = X\) and \(X +_J 1 = 1\). Also, \(X \leq \Omega\) and \(X \leq Y\) imply \(X +_J Y = Y\). Thus \(\n{A} +_J \n{A} = \n{A}\), \(\n{A} +_J \Omega = \Omega\) and \(\n{A} +_J A = A\). To finish showing \(\gen{A} = \set{0, \n{A}, \Omega, A, 1}\) is closed under $+_\emptyset$ and $+_{\set{0}}$, it suffices to perform a few calculations. It is easily checked that  
\begin{alignat*}{2}
	A +_\emptyset A &= 1, 
		&\qquad A +_{\set{0}} A &= A.
\end{alignat*}
For example, \(\set{0,1} \in (A +_\emptyset A)^+\) because \(\set{0,1} = \set{0} \cup_\emptyset \set{1}\), where \(\set{0} \in A^+\) and \(\set{1} \in A^+\). 
\(A +_{\set 0} A = A\) because $+_{\set 0}$ is a lattice join operation.
Finally by \propref{C_0,J...C_N-1,J double suit}, \(C_{0,J}(0) = 0\), \(C_{0,J}(\n A) = C_{0,J}(\Omega) = \Omega\), and \(C_{0,J}(A) = C_{0,J}(1) = 1\).

Now consider \(\gen{B} = \set{0, \n{B}, \Omega, B, 1}\). Since $B$ is flat \(B +_\emptyset B = B\). All the other calculations are the same as for $\gen A$.
Similarly \(\gen{C} = \set{0, \n{C}, \Omega, C, 1}\).

To show \(\gen{A,B} = \set{0, \n{A}, \n{B}, \Omega, B, A, 1}\) observe that 
\begin{alignat*}{2}
	A +_\emptyset B &= 1, 
		&\qquad A +_{\set{0}} B &= A. 
\end{alignat*}
Similarly \(\gen{A,C} = \set{0, \n{A}, \n{C}, \Omega, C, A, 1}\).

To show \(\gen{B,C} = \set{0, \n{A}, \n{B}, \n{C}, \Omega, C, B, A, 1}\), observe that 
\begin{alignat*}{2}
	B +_\emptyset C &= 1, 
		&\qquad B +_{\set{0}} C &= A, \\
	\n{B} +_\emptyset \n{C} &= \Omega,	
		&\qquad \n{B} +_{\set{0}} \n{C} &= \Omega.
\end{alignat*}

Finally, note that if $\D$ is a subalgebra of \(\Cyls_{\mathrm{IFG}_{1}}(\mathbf{2})\) that includes $\norm{v_0 = 0}$, then \(\norm{v_0 = 1} = \n{\norm{v_0 = 0}} \in \D\). Thus $\D$ includes all of the perfect elements in \(\Cyls_{\mathrm{IFG}_{1}}(\mathbf{2})\). Hence \(\D = \Cyls_{\mathrm{IFG}_{1}}(\mathbf{2})\). Similarly, if \(\norm{v_0 = 1} \in \D\), then \(\D = \Cyls_{\mathrm{IFG}_{1}}(\mathbf{2})\).
\end{proof}

\begin{thm} 
\label{CsIFG(2) hereditarily simple}
$\Cyls_{\mathrm{IFG}_{1}}(\mathbf{2})$ is hereditarily simple.
\end{thm}

\begin{proof}
It follows from \lemref{congruence 0, Omega, 1} and \lemref{congruence X < Omega < Y} that the subalgebras \(\set{0, \Omega, 1}\), $\gen{A}$, $\gen{B}$, and $\gen{C}$ are all simple. To show the subalgebra $\gen{A, B}$ is simple, by \lemref{congruence Omega < X < Y} it suffices to show that the congruence $\Cg(A,B)$ generated by $A$ and $B$ is the total congruence. Observe that if \(A \cong B\), then \(1 = A +_\emptyset A \cong B +_\emptyset B = B\), so $\Cg(A,B)$ is the total congruence. A similar argument shows that the subalgebra $\gen{A,C}$ is simple. Finally, to prove the subalgebra $\gen{B,C}$ and $\Cyls_{\mathrm{IFG}_{1}}(\mathbf{2})$ are simple it suffices to show that the congruences $\Cg(A,B)$ and $\Cg(A,C)$ are both the total congruence. But the calculations are the same as before, so we are done.
\end{proof}

\section{$\Cyls_{\mathrm{IFG}_{1}}(\mathbf{3})$ is not hereditarily simple}

Let $\mathbf{3}$ be the structure with universe $\set{0,1,2}$ in which all three elements are named by constant symbols. Then \(\Cyls_{\mathrm{IFG}_{1}}(\mathbf{3}) = \DSuit_1(\set{0,1,2})\), which has 55 elements. Part of the lattice structure of $\Cyls_{\mathrm{IFG}_{1}}(\mathbf{3})$ is shown in Figure \ref{CsIFG(3)}. For simplicity, we only show the interval above $\Omega$. Furthermore, we omit the falsity coordinate and denote each truth coordinate by listing the maximal winning teams. For example, the vertex labeled $\set{0,1}, \set{2}$ denotes the element \(\tuple{\powerset(\set{0,1}) \cup \powerset(\set{2}), \set{\emptyset}}\), and the vertex labeled $\emptyset$ denotes \(\tuple{\set{\emptyset}, \set{\emptyset}} = \Omega\). Readers familiar with the cover of \cite{Balbes:1974} will recognize that Figure \ref{CsIFG(3)} is isomorphic to the free distributive 1-lattice on three generators. To obtain the full lattice structure of \(\Cyls_{\mathrm{IFG}_{1}}(\mathbf 3)\) it is necessary to flip the figure upside-down to get the interval below $\Omega$, then fill in the sides with every possible double suit incomparable to $\Omega$.

\begin{figure} 
\[
\xymatrix{
						&													&	\set{0,1,2} \\
						&													&	\set{0,1}, \set{0,2}, \set{1,2} \ar@{-}[u] \\
						&	\set{0,1}, \set{0,2} \ar@{-}[ur]				&	\set{0,1}, \set{1,2} \ar@{-}[u]										& 	\set{0,2}, \set{1,2} \ar@{-}[ul] \\
						&	\set{0,1}, \set{2} \ar@{-}[u] \ar@{-}[ur]|\hole	& 	\set{0,2}, \set{1} \ar@{-}[ul] \ar@{-}[ur]							& 	\set{1,2}, \set{0} \ar@{-}[ul]|\hole \ar@{-}[u] \\	
\set{0,1} \ar@{-}[ur] 	&													&	\set{0}, \set{1}, \set{2} \ar@{-}[ul] \ar@{-}[u] \ar@{-}[ur]		&	\set{0,2} \ar@{-}[ul]								&	\set{1,2} \ar@{-}[ul] \\
						&	\set{0}, \set{1} \ar@{-}[ul] \ar@{-}[ur] 		& 	\set{0}, \set{2} \ar@{-}[u] \ar@{-}[ur]								& 	\set{1}, \set{2} \ar@{-}[ul] \ar@{-}[ur] \\
						&	\set{0} \ar@{-}[u] \ar@{-}[ur]|\hole 			& 	\set{1} \ar@{-}[ul] \ar@{-}[ur]										& 	\set{2} \ar@{-}[ul]|\hole \ar@{-}[u] \\	
						&													& 	\emptyset\ar@{-}[ul] \ar@{-}[u] \ar@{-}[ur]
}
\]
  \caption{The interval above $\Omega$ in $\Cyls_{\mathrm{IFG}_{1}}(\mathbf{3})$}
  \label{CsIFG(3)}
\end{figure}

The goal of this section is to show that \(\Cyls_{\mathrm{IFG}_{1}}(\mathbf{3})\) is simple, but not hereditarily simple. In fact, every IFG$_N$-algebra whose universe is the collection of all double suits over a set $A$ is simple.

\begin{prop} 
\label{finite constant structures have simple cylindric set algebras}
\(\DSuit(\^NA)\) is simple.
\end{prop}

\begin{proof}
Suppose $X$ and $Y$ are distinct elements of 
\(\DSuit(\^NA)\) such that \(X \cong Y\). Without loss of generality we may assume that there exists a \(V \in Y^+ \setminus X^+\). Let \(Z = \tuple{\powerset(\^NA \setminus V), \powerset(V)}\). Since \(V \notin X^+\) we know that for every \(U \in X^+\) there is an \(\vec a \in V \setminus U\). Hence \(U \cup (\^NA\setminus V) \not= \^NA\). Thus \(1 \not= X +_\emptyset Z \cong Y +_\emptyset Z = 1\). Therefore $\cong$ is the total congruence.
\end{proof}

Recall that in the proof that $\Cyls_{\mathrm{IFG}_{1}}(\mathbf{2})$ is hereditarily simple, we used the fact that \(A +_\emptyset A = 1\) but \(B +_\emptyset B \not= 1\). For any element $X$ of an IFG-algebra, let $nX$ be a abbreviation for \( \underbrace{X +_\emptyset \cdots +_\emptyset X}_{n} \).

\begin{defn} 
The \emph{order}\index{order of an element|mainidx} of an element $X$ is the least positive integer $n$ such that \(nX = 1\). If no such positive integer exists then the order of $X$ is infinite.
\end{defn}

\begin{lem} 
\label{congruent different order}
Let $\cong$ be a congruence on a double-suited IFG$_N$-algebra. If any two elements of different order are congruent, then $\cong$ is the total congruence.
\end{lem}

\begin{proof}
Let \(X \cong Y\). If the order of $X$ is less than the order of $Y$, then for some positive integer $n$ we have \(1 = nX \cong nY \not= 1\).
\end{proof}

We know by \propref{all double suits} and \propref{finite constant structures have simple cylindric set algebras} that $\Cyls_{\mathrm{IFG}_{1}}(\mathbf{3})$ is simple, but we can verify this directly by using the lemmas and the technique of perspective edges. For example, if \(\set{0},\set{1} \cong \set{0,1}\), then \(\set{0}, \set{1}, \set{2} \cong \set{0,1},\set{2}\). But \(\set{0}, \set{1}, \set{2}\) has order 2, while \(\set{0,1},\set{2}\) has order 1, so by \lemref{congruent different order} we have that $\cong$ is the total congruence.

\begin{thm} 
\label{CsIFG_1(3) is not hereditarily simple}
\(\Cyls_{\mathrm{IFG}_{1}}(\mathbf{3})\) is not hereditarily simple.
\end{thm}

\begin{proof}
Let \(A = \tuple{\powerset(\set{0,1}),\ \set{\emptyset}}\) and \(B = \tuple{\powerset(\set{0}) \cup \powerset(\set{1}),\ \set{\emptyset}}\). The subalgebra \(\gen{B} = \set{0, \n{A}, \n{B}, \Omega, B, A, 1}\) is closed under 
$+_\emptyset$ and $+_{\set 0}$ because\begin{alignat*}{2}
A +_\emptyset A &= A, 		&\qquad	 	A +_{\set{0}} A &= A, \\
A +_\emptyset B &= A,		&\qquad		A +_{\set{0}} B &= A, \\
B +_\emptyset B &= A,		&\qquad		B +_{\set{0}} B &= B, \\
A +_\emptyset \n{A} &= A,		&\qquad		A +_{\set{0}} \n{A} &= A, \\
A +_\emptyset \n{B} &= A,		&\qquad		A +_{\set{0}} \n{B} &= A, \\
B +_\emptyset \n{A} &= B,		&\qquad		B +_{\set{0}} \n{A} &= B, \\
B +_\emptyset \n{B} &= B,		&\qquad		B +_{\set{0}} \n{B} &= B, \\
\n{A} +_\emptyset \n{A} &= \n{A},	&\qquad		\n{A} +_{\set{0}} \n{A} &= \n{A}, \\
\n{A} +_\emptyset \n{B} &= \n{B},	&\qquad		\n{A} +_{\set{0}} \n{B} &= \n{B}, \\
\n{B} +_\emptyset \n{B} &= \n{B},	&\qquad		\n{B} +_{\set{0}} \n{B} &= \n{B}. 
\end{alignat*}
All of the $+_{\set 0}$ calculations are easy to check by looking at the lattice. The $+_\emptyset$ calculations require some computation. First, \(A +_\emptyset A = A\) because $A$ is flat. Second, \(A +_\emptyset B = A\) because \(\set{0, 1} = \set{0} \cup_\emptyset \set{1}\), where \(\set{0} \in A^+\) and \(\set{1} \in B^+\), while \((A +_\emptyset B)^- = A^- \cap B^- = A^-\). Similarly, \(B +_\emptyset B = A\). The remaining $+_\emptyset$ calculations all follow from \propref{X < Omega,Y} or \propref{flat absorption}. Finally, the set is closed under $C_{0,J}$ by \propref{C_0,J...C_N-1,J double suit}.

Let $\cong$ denote the equivalence relation that makes \(A \cong B\) and \(\n{A} \cong \n{B}\), but makes no other pair of distinct elements equivalent. To verify that $\cong$ is a congruence, observe that $\cong$ is preserved under $\n{\,}$ because \(\n{A} \cong \n{B}\) and \(\n{(\n{A})} = A \cong B = \n{(\n{B})}\). It is preserved under $C_{0,J}$ because \(C_{0,J}(A) = 1 = C_{0,J}(B)\) and \(C_{0,J}(\n{A}) = \Omega = C_{0,J}(\n{B})\). Finally, the calculations above show that $\cong$ is preserved under $+_\emptyset$ and $+_{\set{0}}$.
Thus $\cong$ is a nontrivial, non-total congruence. Therefore $\gen{B}$ is not simple.
\end{proof}


\section{``Iff'' is not expressible in IFG logic}

Let \(\phi \himplies{J} \psi\)\index{$\phi \himplies{J} \psi$} be an abbreviation for \({\hneg\phi} \hor{J} \psi\), and let \(\phi \hiff{J} \psi\)\index{$\phi \hiff{J} \psi$} be an abbreviation for 
\[
(\phi \himplies{J} \psi) \hand{J} (\psi \himplies{J} \phi).
\]
It will be useful to know when \(\A \modeltf \phi \himplies{J} \psi[V]\) and \(\A \modeltf \phi \hiff{J} \psi[V]\). It follows immediately from the definitions that for \(\phi \himplies{J} \psi\), 
\begin{itemize}
	\item[(+)] \(\A \modelt \phi \himplies{J} \psi[V]\) if and only \(\A \modelf \phi[V_1]\) and \(\A \modelt \psi[V_2]\) for some \(V = V_1 \cup_J V_2\),
	\item[($-$)] \(\A \modelf \phi \himplies{J} \psi[W]\) if and only if \(\A \modelt \phi[W]\) and \(\A \modelf \psi[W]\).
\end{itemize}
Similarly for \(\phi \hiff{J} \psi\), 
\begin{itemize}
	\item[(+)] \(\A \modelt \phi \hiff{J} \psi[V]\) if and only if \(\A \modelf \phi[V_1]\) and \(\A \modelt \psi[V_2]\) for some \(V = V_1 \cup_J V_2\), and \(\A \modelt \phi[V_3]\) and \(\A \modelf \psi[V_4]\) for some \(V = V_3 \cup_J V_4\),
	\item[($-$)] \(\A \modelf \phi \hiff{J} \psi[W]\) if and only if \(\A \modelt \phi[W_1]\), \(\A \modelf \psi[W_1]\), \(\A \modelf \phi[W_2]\), and \(\A \modelt \psi[W_2]\) for some \(W = W_1 \cup_J W_2\).
\end{itemize}
In particular, the semantics for \(\phi \himplies{N} \psi\) are
\begin{itemize}
	\item[(+)] \(\A \modelt \phi \himplies{N} \psi[V]\) if and only if \(\A \modelf \phi[V]\) or \(\A \modelt \psi[V]\),
	\item[($-$)] \(\A \modelf \phi \himplies{N} \psi[W]\) if and only if \(\A \modelt \phi[W]\) and \(\A \modelf \psi[W]\),
\end{itemize}
and for \(\phi \hiff{N} \psi\),
\begin{itemize}
	\item[(+)] \(\A \modelt \phi \hiff{N} \psi[V]\) if and only if \(\A \modelt \phi[V]\) and \(\A \modelt\psi[V]\), or \(\A \modelf \phi[V]\) and \(\A \modelf \psi[V]\),
	\item[($-$)] \(\A \modelf \phi \hiff{N} \psi[W]\) if and only if \(\A \modelt \phi[W]\) and 
	\(\A \modelf \psi[W]\), or \(\A \modelf \phi[W]\) and \(\A \modelt \psi[W]\).
\end{itemize}
For example, \(\A \modelt \phi \hiff{N} \psi[V]\) 
if and only if \(\A \modelt ({\hneg\phi} \hor{N} \psi) \hand{N} (\phi \hor{N} {\hneg\psi})[V]\) 
if and only if \(\A \modelt ({\hneg\phi} \hor{N} \psi)[V]\) and \(\A \modelt (\phi \hor{N} {\hneg\psi})[V]\)
if and only if \(\A \modelf \phi[V]\) or \(\A \modelt \psi[V]\), and \(\A \modelt \phi[V]\) or \(\A \modelf \psi[V]\) 
if and only if \(\A \modelf \phi[V]\) and \(\A \modelf \psi[V]\), or \(\A \modelt \psi[V]\) and \(\A \modelt \phi[V]\).

\begin{prop} 
\label{hiff_0}
For any IFG$_N$-formulas $\phi$ and $\psi$, \(\A \modelt \phi \hiff{\emptyset} \psi\) if and only if \(\norm{\phi}_\A = \norm{\psi}_\A\) and both are perfect.
\end{prop}

\begin{proof}
Suppose \(\A \modelt (\phi \hiff{\emptyset} \psi)[\^NA]\). Then there exist \(V, V' \subseteq \^NA\) such that 
\(\A \modelt \phi[V]\), \(\A \modelf \psi[\^NA\setminus V]\), \(\A \modelf \phi[V']\), and \(\A \modelt \psi[\^NA \setminus V']\). Thus \(V \cap V' = \emptyset\) and \((\^NA \setminus V) \cap (\^NA \setminus V') = \emptyset\). Therefore \(V' = \^NA \setminus V\), and \(\norm{\phi}_\A = \tuple{\powerset(V), \powerset(\^NA \setminus V)} = \norm{\psi}_\A\).
\end{proof}

\begin{prop} 
\label{hiff_N}
For any IFG$_N$-formulas $\phi$ and $\psi$, \(\A \modelt \phi \hiff{N} \psi\) if and only if \(\norm{\phi}_\A = \norm{\psi}_\A \in \set{0,1}\).
\end{prop}

\begin{proof}
Suppose \(\A \modelt (\phi \hiff{N} \psi)[\^NA]\). Then \(\A \modelt \phi[\^NA]\) and \(\A \modelt\psi[\^NA]\), in which case \(\norm{\phi}_\A = \norm{\psi}_\A = 1 \), or \(\A \modelf \phi[\^NA]\) and \(\A \modelf \psi[\^NA]\), in which case \(\norm{\phi}_\A = \norm{\psi}_\A = 0\).
\end{proof}

\begin{defn} 
An \emph{IFG$_N$-schema}\index{IFG$_N$-schema|mainidx} involving $k$ formula variables is an element of the smallest set $\Xi$ satisfying the following conditions. 
\begin{enumerate}
	\item The formula variables \(\alpha_0, \ldots, \alpha_{k-1}\) belong to $\Xi$.
	\item For all \(i,j < N\), the formula \(v_i = v_j\) belongs to $\Xi$.
	\item If $\xi$ belongs to $\Xi$, then ${\hneg\xi}$ belongs to $\Xi$.
	\item If $\xi_1$ and $\xi_2$ belong to $\Xi$, and \(J \subseteq N\), then \(\xi_1 \hor{J} \xi_2\) belongs to $\Xi$.
	\item If $\xi$ belongs to $\Xi$, \(n < N\), and \(J \subseteq N\), then \(\hexists{v_n}{J}\xi\) belongs to $\Xi$.
\end{enumerate}
Note that the symbols \(\alpha_0, \ldots, \alpha_{k-1}\) are distinct from the usual variables $v_0$, \ldots, $v_{k-1}$.  If $\xi$ is an IFG$_N$-schema involving $k$ formula variables, and \(\phi_0, \ldots, \phi_{k-1}\) are IFG$_N$-formulas, then the IFG$_N$-formula \(\xi(\phi_0, \ldots, \phi_{k-1})\) is called an \emph{instance}\index{instance of a schema|mainidx} of $\xi$.
\end{defn}

\begin{defn} 
Every IFG$_N$-schema $\xi$ has a corresponding term $T_\xi$ in the language of IFG$_N$-algebras. $T_\xi$ is defined recursively as follows:
\begin{enumerate}
	\item \(T_{\alpha_i} = X_i\),
	\item \(T_{v_i = v_j} = D_{ij}\),
	\item \(T_{\hneg\,\xi} = \n{(T_\xi)}\),
	\item \(T_{\xi_1 \hor{J} \xi_2} = T_{\xi_1} +_J T_{\xi_2}\),
	\item \(T_{\hexists{v_n}{J}\xi} = C_{n,J}(T_\xi)\).
\end{enumerate}
\end{defn}

\begin{lem} 
\label{schema terms}
Let $\xi$ be an IFG$_N$-schema involving $k$ formula variables, and let $T_\xi$ be its corresponding term. Then for any IFG$_N$-formulas \(\phi_0, \ldots, \phi_{k-1}\) and any suitable structure $\A$,
\[
\norm{\xi(\phi_0, \ldots, \phi_{k-1})} = T_\xi^{\Cyls_{\mathrm{IFG}_{N}}(\A)}(\norm{\phi_0}, \ldots, \norm{\phi_{k-1}}).
\]
\end{lem}

\begin{proof}
If $\xi$ is a formula variable $\alpha_i$, then \(T_\xi = X_i\), so
\[
\norm{\xi(\phi_0, \ldots, \phi_{k-1})} = \norm{\phi_i} = T_\xi^{\Cyls_{\mathrm{IFG}_{N}}(\A)}(\norm{\phi_0}, \ldots, \norm{\phi_{k-1}}).
\]
If $\xi$ is \(v_i = v_j\), then \(T_\xi = D_{ij}\), so 
\[
\norm{\xi(\phi_0, \ldots, \phi_{k-1})} = \norm{v_i = v_j} = D_{ij}^{\Cyls_{\mathrm{IFG}_{N}}(\A)} = T_\xi^{\Cyls_{\mathrm{IFG}_{N}}(\A)}(\norm{\phi_0}, \ldots, \norm{\phi_{k-1}}).
\]
Now assume that 
{\allowdisplaybreaks
\begin{align*}
	\norm{\xi(\phi_0, \ldots, \phi_{k-1})} 
		&= T_\xi^{\Cyls_{\mathrm{IFG}_{N}}(\A)}(\norm{\phi_0}, \ldots, \norm{\phi_{k-1}}), \\
	\norm{\xi_1(\phi_0, \ldots, \phi_{k-1})}
		&= T_{\xi_1}^{\Cyls_{\mathrm{IFG}_{N}}(\A)}(\norm{\phi_0}, \ldots, \norm{\phi_{k-1}}), \\
	\norm{\xi_2(\phi_0, \ldots, \phi_{k-1})} 
	&= T_{\xi_2}^{\Cyls_{\mathrm{IFG}_{N}}(\A)}(\norm{\phi_0}, \ldots, \norm{\phi_{k-1}}).
\intertext{Then} 
	\norm{\hneg\xi(\phi_0, \ldots, \phi_{k-1})} 
		&= \n{\norm{\xi(\phi_0, \ldots, \phi_{k-1})}} \\
		&= \n{\left(T_\xi^{\Cyls_{\mathrm{IFG}_{N}}(\A)}(\norm{\phi_0}, \ldots, \norm{\phi_{k-1}})\right)} \\
		&= T_{\hneg\,\xi}^{\Cyls_{\mathrm{IFG}_{N}}(\A)}(\norm{\phi_0}, \ldots, \norm{\phi_{k-1}}), \\[10pt]
	\norm{\xi_1 \hor{J} \xi_2(\phi_0, \ldots, \phi_{k-1})} 
		&= \norm{\xi_1(\phi_0, \ldots, \phi_{k-1})} +_J \norm{\xi_2(\phi_0, \ldots, \phi_{k-1})} \\
		&= T_{\xi_1}^{\Cyls_{\mathrm{IFG}_{N}}(\A)}(\norm{\phi_0}, \ldots, \norm{\phi_{k-1}}) \\
		&\qquad +_J T_{\xi_2}^{\Cyls_{\mathrm{IFG}_{N}}(\A)}(\norm{\phi_0}, \ldots, \norm{\phi_{k-1}}) \\
		&= T_{\xi_1 \hor{J} \xi_2}^{\Cyls_{\mathrm{IFG}_{N}}(\A)}(\norm{\phi_0}, \ldots, \norm{\phi_{k-1}}), \\[10pt]
	\norm{\hexists{v_n}{J}\xi(\phi_0, \ldots, \phi_{k-1})} 
		&= C_{n,J}\norm{\xi(\phi_0, \ldots, \phi_{k-1})} \\
		&= C_{n,J}\left(T_\xi^{\Cyls_{\mathrm{IFG}_{N}}(\A)}(\norm{\phi_0}, \ldots, \norm{\phi_{k-1}})\right) \\
		&= T_{\hexists{v_n}{J}\xi}^{\Cyls_{\mathrm{IFG}_{N}}(\A)}(\norm{\phi_0}, \ldots, \norm{\phi_{k-1}}). 
\end{align*}}
\end{proof}

\begin{prop} 
\label{term operation implies hereditarily simple}
Any double-suited IFG-algebra that has a term operation $T(X,Y)$ such that 
\(T(X,Y) = 1\) if and only if \(X = Y\) is hereditarily simple.
\end{prop}

\begin{proof}
Suppose $\C$ is a double-suited IFG-algebra that has such a term oper\-ation. Then for any \(X \not= Y\) we have \(\tuple{1, Z} = \tuple{T(X,X), T(X,Y)} \in \Cg(X,Y)\), where $Z$ is some element different than 1. Hence $\Cg(X,Y)$ is the total congruence. Thus $\C$ is simple. 
Furthermore, the sentence 
\[
\forall X \forall Y[ T(X,Y) = 1 \iff X = Y]
\]
is universal, and so must hold in every subalgebra of $\C$. Hence $\C$ is hereditarily simple.
\end{proof}

\begin{thm} 
\label{no schema}
There is no IFG$_1$-schema $\xi$ involving two formula variables such that for every pair of IFG$_1$-formulas $\phi$ and $\psi$, and every suitable structure $\A$, we have
\[
\A \modelt \xi(\phi, \psi) \quad \text{iff} \quad \norm{\phi}_\A = \norm{\psi}_\A.
\]
\end{thm}

\begin{proof}
Suppose $\xi$ were such a schema. Then the corresponding term $T_\xi$ would have the property that for any $\A$ and any \(\norm{\phi}_\A, \norm{\psi}_\A \in \Cyls_{\mathrm{IFG}_{1}}(\A)\), 
\begin{align*}
	T_\xi^{\Cyls_{\mathrm{IFG}_{1}}(\A)}(\norm{\phi}_\A, \norm{\psi}_\A) = 1 \quad 
		&\text{iff} \quad \norm{\xi(\phi,\psi)}_\A = 1 \\
		&\text{iff} \quad \A \modelt \xi(\phi,\psi) \\
		&\text{iff} \quad \norm{\phi}_\A = \norm{\psi}_\A.
\end{align*}
Thus every \(\Cyls_{\mathrm{IFG}_{1}}(\A)\) would be hereditarily simple. However \(\Cyls_{\mathrm{IFG}_{1}}(\mathbf{3})\) is not hereditarily simple.
\end{proof}

	

\end{document}